\documentclass[]{article}
\begin{document}
\newtheorem{proposition}{Proposition}[section]
\newtheorem{definition}{Definition}[section]
\newtheorem{lemma}{Lemma}[section]

\title{\bf How to Define Automatic Differentiation  }
\author{Keqin Liu \\Department of Mathematics\\The University of British Columbia\\Vancouver, BC\\
Canada, V6T 1Z2}
\date{June, 2020}
\maketitle

\begin{abstract} Based on a class of associative algebras with zero-divisors which are called real-like algebras by us, we  introduce a way of defining automatic differentiation  and present different ways of doing automatic differentiation to compute  the first, the second and the third derivatives of a function  exactly and simultaneously.

\end{abstract}

\bigskip
Automatic differentiation is indispensable to perform many optimization algorithms for 
large-scale machine learning problems. Computer scientists in machine learning community describe 
automatic differentiation as the computing techniques to evaluate accurately the numerical values of a derivative up to machine precision without using the explicit formula for the derivative, and emphasize that automatic differentiation is neither numerical differentiation nor symbolic differentiation. These popular description and emphasis about automatic differentiation are informative, but they do not indicate which kinds of mathematical objects we should search to find the techniques of doing automatic differentiation. Hence, an interesting problem is how to conceptualize automatic differentiation mathematically so that we know exactly which kinds of mathematical objects are sufficient for getting a procedure of executing automatic differentiation. 

\medskip
In addition to machine learning, many applications in scientific computation require to  compute higher-order derivatives precisely. It is well-known that even a differential function is given explicitly by a mathematical formula, the explicit expression of its higher-order derivatives is usually too complicated to evaluate the higher-order derivatives precisely. Hence, in order to compute higher-order derivatives at machine precision, automatic differentiation seems to be  unavoidable. The success of first-order automatic differentiation in scientific computation comes from its nice property  given in the section 3.1.1 of \cite{GPRS}:  the derivative is obtained exactly by computing the value of a function from the dual  numbers to the dual numbers  just once, where the dual numbers were introduced by C. L. Clifford in \cite{C}.
To the best of our knowledge, there is only one way of doing  higher order automatic differentiation in  scientific computation.
Therefore, another interesting problem  is to find out if there exists a different way of doing higher order automatic differentiation so that the  higher-order derivatives up to an arbitrary order  can be obtained exactly by computing the value of a function from an algebra to the algebra just once.

\medskip
The goal of this paper is to give our answer to these two problems by using a class of associative algebras called real-like  algebras. Based on our observation, the automatic differentiation techniques appearing in machine learning community are always related to some  real-like  algebras which generalize the dual numbers.
In section 1 of this paper, we define real-like algebras and discuss their basic properties.  
In section 2, we introduce the concept of the $n$-th order automatic differentiation induced by a real-like algebra and 
explain why more work should be done to have a solid mathematical foundation of higher order automatic differentiation used in scientific computation.
To simplify our presentation, we just discuss  the $3$-rd order automatic differentiation in this paper.
After  proving the proposition which presents many ways of doing the $3$-rd order automatic differentiation induced by a real-like algebra in the last section of this paper,  we do an example to demonstrate how to use  two different ways of implementing  the $3$-rd order automatic differentiation technique to get the first, the second and the third derivatives of a function  exactly and simultaneously.

\bigskip
Throughout this paper,   the real number field is denoted by  $\mathcal{R}$ , the range of a function $f$ is denoted by $Im f$, and an associative algebra with the identity is just called  a unital  associative algebra.

\bigskip
\section{Real-like Algebras}

Let $n$ be non-negative integer. We say that a vector space $A$ over the real number field 
$\mathcal{R}$  is the {\bf direct sum} of its non-zero subspaces $A_0$, $A_1$, \dots  , $A_n$ and we write $A=\displaystyle\bigoplus_{i=0}^n A_i$ if 
each $x$ in A can be represented uniquely in the form $x=\displaystyle\sum _{i=0}^n x_i$ for $x_i\in A_i$ for $0\le i\le n$. The subspaces $A_0$, $A_1$, \dots  , $A_n$ are called the {\bf 
homogeneous subspaces} of $A$. The elements of $A_i$ are said to be {\bf homogeneous of degree} $i$ for $0\le i\le n$. After expressing an element $a$ in $A$ as a sum of non-zero homogeneous elements of distinct degrees, these non-zero homogeneous elements are called the {\bf homogeneous components} of $a$ and the  homogeneous components of $a$ of least degree is called the {\bf initial component} of $a$.

\medskip
We now define real-like algebras which are the kind of real associative algebras we need in the study of automatic differentiation. 

\medskip
\begin{definition}\label{def1.1} A commutative associative algebra $A$  over the real number field $\mathcal{R}$ is called a {\bf real-like algebra} if $A$ is the direct sum $A=\displaystyle\bigoplus_{i=0}^n A_i$ of its non-zero subspaces $A_0$, $A_1$, \dots  , $A_n$ satisfying $A_0=\mathcal{R}$ and 
$ A_iA_j\subseteq A_{i+j}$ for $0\le i, \, j\le n$, where $A_{i+j}:=0$ for $i+j>n$.
\end{definition}

\bigskip
Based on our observation, the automatic differentiation techniques appearing in scientific computation are always related to some  real-like  algebras which generalize the dual numbers. Let us give two examples of real-like  algebras.

\medskip
{\bf Example 1} One example  is  the  $(n+1)$-dimensional truncated polynomial real algebra  $\mathcal{R}^{(n+1)}$, where 
$$\mathcal{R}^{(n+1)}:=\mathcal{R}[X]/<\{X^k\,|\, k> n\}>$$ 
is the quotient associative algebra of the polynomial ring  $\mathcal{R}[X]$ with respect to the ideal 
$<\{X^k\,|\, k>n\}>$ generated by the subset $\{X^k\,|\, k> n\}$ of $\mathcal{R}[X]$. Clearly, $\mathcal{R}^{(n+1)}$ is a real-like algebra. In fact, if 
$\varepsilon:=X+<\{X^k\,|\, k> n\}>$,
then we have $\mathcal{R}^{(n+1)}=\displaystyle\bigoplus_{i=0}^{n}A_i$ and
$$
A_i:= \mathcal{R}\varepsilon ^i, \quad\mathcal{R}\varepsilon ^0=\mathcal{R}, \quad
(\mathcal{R}\varepsilon ^i)(\mathcal{R}\varepsilon ^j)
=\left\{\begin{array}{cl}\mathcal{R}\varepsilon ^{i+j}&
\mbox{if $i+j\le n$}\\0&\mbox{if $i+j> n$}\end{array}\right.,
$$
where $0\le i, \, j\le n$. 

\bigskip
{\bf Example 2} In addition to the  $(n+1)$-dimensional  truncated polynomial real algebras,  another example of  the real-like algebras is the 4-dimensional associative algebra  
${\bf A}=\displaystyle\bigoplus_{i=0}^3 A_i$, where 
$A_0=\mathcal{R}$,  $A_1=\mathcal{R}\epsilon _1$, $A_2=\mathcal{R}\epsilon _2$, 
$A_3=\mathcal{R}\epsilon _1\epsilon _2$
and the multiplication  is given by $\epsilon _1^2=\epsilon _2^2=0$ and
$\epsilon _1\epsilon _2=\epsilon _2\epsilon _1$.  This 4-dimensional real-like algebra
${\bf A}$ has been used for exact  second derivative calculations in \cite{FA}.
 
\bigskip
The  $(n+1)$-dimensional truncated polynomial real algebra  $\mathcal{R}^{(n+1)}$, which is denoted by 
by $_nD_1$ in  Section 2.2.2 of \cite{B} and  by $\mathcal{P}_{n+1}$ in  Section 13.2 of \cite{GW} respectively,  has appeared in  automatic differentiation for a long time.   To the best of our knowledge, although $\mathcal{R}^{(n+1)}$ is a $\ell_1$-normed algebra by \cite{B} or by \cite{GW}, a $\ell_2$-normed algebraic structure has not been introduced on the truncated polynomial real algebra  $\mathcal{R}^{(n+1)}$. Since $\ell_2$-norm is generally preferred in neural networks and more computational efficient than $\ell_1$-norm, it is advantageous to have a 
$\ell_2$-normed algebraic structure on the algebras appearing in the study of automatic differentiation. At the end of  this section,  we give  many ways of introducing a 
$\ell_2$-norm which make  the truncated polynomial real algebra  $\mathcal{R}^{(n+1)}$ into a normed algebra.
For convenience, we use {\bf $\mathcal{R}^{(n+1)}$-numbers} to name the elements of the truncated polynomial real algebra  $\mathcal{R}^{(n+1)}$. Clearly, $\mathcal{R}^{(2)}$-number are  the dual numbers introduced by C. L. Clifford in \cite{C}.  Based on our research about the applications of the dual numbers, we strongly feel that if there exists a class of new numbers which can be used to extend the known mathematics based on real numbers in a satisfyingly way, then $\mathcal{R}^{(n+1)}$-numbers should be the best candidate for this class of new numbers. 
 
\medskip
The following proposition gives the basic properties of real-like algebras. 

\begin{proposition}\label{pr1.1} Let $A=\displaystyle\bigoplus_{i=0}^n A_i$ be a real-like algebra and let $a=\displaystyle\sum_{i=0}^n a_i$ be an element of $A$ with $a_i\in A_i$ for $0\le i\le n$.
\begin{description} 
\item[(i)] $a$ is a zero-divisor  if and only if $a_0=0$.
\item[(ii)] $a$ is invertible  if and only if $a_0\ne 0$, in which case, the inverse $a^{-1}$ of $a$ is given by
$a^{-1}=\displaystyle\sum_{i=0}^n \displaystyle\frac{det\, M_i}{a_0^{n+1}}$, where $M_i$ is the  $(n+1)\times (n+1)$-matrix obtained by replaying the $i$-th column of the  $(n+1)\times (n+1)$-matrix 
$$
M:=\left[\begin{array}{cccccc}
a_0&0&0&\cdots&0&0\\
a_1&a_0&0&\cdots&0&0\\
a_2&a_1&a_0&\cdots&0&0\\
\vdots&\vdots&\vdots&\cdots&\vdots&\vdots\\
a_{n-1}&a_{n-2}&a_{n-3}&\cdots&a_0&0\\
a_n&a_{n-1}&a_{n-2}&\cdots&a_1&a_0
\end{array}\right]
$$  
with the $(n+1)\times 1$-matrix $\left[\begin{array}{c}
1\\0\\\vdots\\0
\end{array}\right]$ and $det\, M_i$ is the determinant of the $(n+1)\times (n+1)$-matrix $M_i$.
\end{description}
\end{proposition}

\medskip
\noindent
{\bf Proof} (i) If $a$ is a zero-divisor, then $ab=0$ for some $0\ne b\in A$. Let $p$ be the degree of the initial component of $b$. Then we have $b=\displaystyle\sum_{i=p}^n b_i$, where $b_i\in A_i$ for $p\le i\le n$ and $b_p\ne 0$. Assume that $a_0\ne 0$. By the fact that $a_0$ is in $A_0=\mathcal{R}$,  the inverse $a_0^{-1}$ of $a_0$ exists. It follows that
\begin{equation}\label{eq1}
0=a_0^{-1}ab=a_0^{-1}\Big(\displaystyle\sum_{i=0}^n a_i\Big)(\displaystyle\sum_{i=p}^n b_i)
=b_p+\underbrace{\displaystyle\sum_{i=p+1}^n b_i+
\Big(\displaystyle\sum_{i=1}^n a_i\Big)(\displaystyle\sum_{i=p}^n b_i)}.
\end{equation}
Since the  degree of the initial component of 
$\displaystyle\sum_{i=p+1}^n b_i+
\Big(\displaystyle\sum_{i=1}^n a_i\Big)(\displaystyle\sum_{i=p}^n b_i)$
is at least $p+1$, we have to have $b_p=0$ by (\ref{eq1}), which is impossible. This proves that 
$a_0$ has to be $0$.

\medskip
Conversely, if $a_0= 0$, then $a=\displaystyle\sum_{i=1}^n a_i$. After choosing 
$0\ne b_n\in A_n$, we get $ab_n=\Big(\displaystyle\sum_{i=1}^n a_i\Big)b_n=0$. This proves that 
$a$ is a zero-divisor.

\bigskip
(ii) $a$ is invertible if and only if there exists  $x=\displaystyle\sum_{i=0}^n x_i$ with $x_i\in A_i$ for $0\le i\le n$ such that $\Big(\displaystyle\sum_{i=0}^n a_i\Big)\Big(\displaystyle\sum_{i=0}^n x_i\Big)=ax=1$, which is equivalent to
$$
1=a_0x_0+(a_0x_1+a_1x_0)+(a_2x_0+a_1x_1+a_0x_2)+\cdots +(a_nx_0+a_{n-1}x_1+\cdots +a_0x_n)
$$
or
$$
M\,\left[\begin{array}{c}
x_0\\x_1\\x_2\\\vdots\\x_{n-1}\\x_n
\end{array}\right]=
\left[\begin{array}{cccccc}
a_0&0&0&\cdots&0&0\\
a_1&a_0&0&\cdots&0&0\\
a_2&a_1&a_0&\cdots&0&0\\
\vdots&\vdots&\vdots&\cdots&\vdots&\vdots\\
a_{n-1}&a_{n-2}&a_{n-3}&\cdots&a_0&0\\
a_n&a_{n-1}&a_{n-2}&\cdots&a_1&a_0
\end{array}\right]
\left[\begin{array}{c}
x_0\\x_1\\x_2\\\vdots\\x_{n-1}\\x_n
\end{array}\right]
=\left[\begin{array}{c}
1\\0\\0\\\vdots\\0\\0
\end{array}\right].
$$
It follows from that $x_i=\displaystyle\frac{det\, M_i}{a_0^{n+1}}$ for $0\le i\le n$. Thus (ii) holds.

\hfill\raisebox{1mm}{\framebox[2mm]{}}

\bigskip
To study the norm algebraic structure of a real-like algebra, we introduce the concept of a homogeneous norm  in the following

\medskip
\begin{definition}\label{def1.2} We say that a real-like algebra $A=\displaystyle\bigoplus_{i=0}^n A_i$ has a {\bf homogeneous norm $|\,\,|_*$} if $|\,\,|_*: \displaystyle\bigcup_{i=0}^n A_i\to \mathcal{R}$ is a function on the set 
$\displaystyle\bigcup_{i=0}^n A_i$ of homogeneous elements of $A$ such that  for $r\in\mathcal{R}$, $x_i\in A_i$,  $y_j\in A_j$ and $0\le i,\, j\le n$: \begin{description} 
\item[(i)] $|x_i|_*\ge 0$, and $|x_i|_*= 0$ if and only if $x_i= 0$,
\item[(ii)] $|rx_i|_*=|r|\,|x_i|_*$,
\item[(iii)] $|x_i+y_i|_*\le |x_i|_*+|y_i|_*$,
\item[(iv)] $|x_iy_j|_*\le |x_i|_*\,|y_j|_*$.
\end{description}
\end{definition}

\medskip
Let $A=\displaystyle\bigoplus_{i=0}^n A_i$ be a  real-like algebra which has a homogeneous norm $|\,\,|_*: \displaystyle\bigcup_{i=0}^n A_i\to \mathcal{R}$. Mimicking the definitions of the ordinary $\ell_1$-norm and $\ell_2$-norm on the $(n+1)$-dimensional Euclidean space $\mathcal{R}^{n+1}$, we have the  following natural extensions 
$||\,\,||_1: A\to \mathcal{R}$ and $||\,\,||_{2,*}: A\to \mathcal{R}$ of the function 
$|\,\,|_{*}: \displaystyle\bigcup_{i=0}^n A_i\to \mathcal{R}$:
\begin{equation}\label{eq2}
||a||_1:=\displaystyle\sum_{i=0}^n |a_i|_*\quad\mbox{and}\quad 
||a||_{2,*}:=\Big(\displaystyle\sum_{i=0}^n |a_i|_*^2\Big)^{\frac 12},
\end{equation}
where  $a=\displaystyle\sum_{i=0}^n a_i\in A$ with $a_i\in A_i$ for $0\le i\le n$. The following proposition gives the basic properties of the two real-valued functions $||\,\,||_1$ and $||\,\,||_{2,*}$.

\medskip
\begin{proposition}\label{pr1.2} Let $A=\displaystyle\bigoplus_{i=0}^n A_i$ be a real-like algebra which has a homogeneous norm $|\,\,|_*: \displaystyle\bigcup_{i=0}^n A_i\to \mathcal{R}$, and let $||\,\,||_1$ and $||\,\,||_{2,*}$ be the real-valued functions defined by (\ref{eq2}).
\begin{description} 
\item[(i)]  $A$ is a normed algebra with respect to the norm $||\,\,||_1$. 
\item[(ii)] If $n\ge 1$ and $|e|_*=1$ for the identity $e$ of the algebra $A$, then $A$ can not be made into  a normed algebra via the real-valued function $||\,\,||_{2,*}$. 
\item[(iii)] $||\,\,||_{2,*}$ is a norm on $A$.
\end{description}
\end{proposition}

\medskip
\noindent
{\bf Proof}  Recall that a function $||\,\,||: A\to \mathcal{R}$ is called a {\bf norm} on $A$ if for $r\in\mathcal{R}$ and $x$, $y\in A$, we have
\begin{equation}\label{eq4}
||rx||=|r|\,||x||, \quad ||x||\ge 0, \quad\mbox{and $||x||= 0$ if and only if $x= 0$}
\end{equation}
and
\begin{equation}\label{eq5}
||x+y||\le ||x||+||y||.
\end{equation}
Also, $A$ is called a {\bf normed algebra} if there is a norm $||\,\,||$ on $A$ such that
\begin{equation}\label{eq6}
||xy||\le ||x||\,||y||\quad\mbox{for all $x$, $y\in A$.}
\end{equation}

\bigskip
The proof of Proposition \ref{pr1.2} follows from direct computations. We now proof 
(ii) to explain the way of doing the computation. Since $n\ge 1$, we have $a_1\in A_1$ such that  $ a_1\ne 0$. Then 
\begin{equation}\label{eq7}
|a_1|_*> 0\quad\mbox{and}\quad |a_1^2|_*=|a_1\cdot a_1|_*\le |a_1|_*\cdot |a_1|_*=|a_1|_*^2.
\end{equation}
It follows from (\ref{eq2}) and (\ref{eq7}) that
\begin{eqnarray*}
&&||(e+a_1)\cdot (e+a_1)||_{2,*}^2=||e+2a_1+a_1^2||_{2,*}^2=|e|_*^2+|2a_1|_*^2+|a_1^2|_*^2\\
&=&|e|_*^2+(2\,|a_1|_*)^2+|a_1^2|_*^2=|e|_*^2+4\,|a_1|_*^2+|a_1^2|_*^2
>|e|_*^2+2\,|a_1|_*^2+|a_1^2|_*^2\\
&\ge & |e|_*^2+2\,|a_1^2|_*+|a_1^2|_*^2 = \big(|e|_*+|a_1^2|_*\big)^2
=\big(|e|_*^2+|a_1^2|_*\big)^2=\big(||e+a_1||_{2,*}^2\big)^2\\
&=&
\big(||e+a_1||_{2,*}\cdot ||e+a_1||_{2,*}\big)^2
\end{eqnarray*}
or 
$$||(e+a_1)\cdot (e+a_1)||_{2,*}>||e+a_1||_{2,*}\cdot ||e+a_1||_{2,*},$$
which proves that (\ref{eq6}) fails for $x=y=e+a_1$.

\hfill\raisebox{1mm}{\framebox[2mm]{}}

\bigskip
Obviously,  the map 
$|\,\,|_*: \displaystyle\bigcup_{i=0}^{n}\mathcal{R}\varepsilon ^i\to \mathcal{R}$ defined by 
$$|x\varepsilon ^i|_*: =|x|\quad\mbox{for $x\in \mathcal{R}$ and $0\le i\le n$}$$ 
is a  homogeneous norm on the truncated polynomial real algebra $\mathcal{R}^{(n+1)}$ which satisfies the assumption in Proposition \ref{pr1.2} (ii), where $|x|$ is the absolute value of the real number $x$. Hence, the natural idea of extending the ordinary way  of defining a $\ell_2$-norm on  the $(n+1)$-dimensional Euclidean space 
$\mathcal{R}^{n+1}$ can not give a $\ell_2$-normed algebraic structure on the truncated polynomial real algebra $\mathcal{R}^{(n+1)}$  by Proposition \ref{pr1.2} (ii). This is possibly why we have not seen the way of making the truncated polynomial real algebra $\mathcal{R}^{(n+1)}$ into  a $\ell_2$-normed algebra in  automatic differentiation community even it  has a $\ell_2$-normed algebraic structure. We now give many ways of introducing a $\ell_2$-normed algebraic structure on the truncated polynomial real algebra 
$\mathcal{R}^{(n+1)}$.

\medskip
\begin{proposition}\label{pr1.3} Let $\beta$ be a positive  real number. If  $||\,\,||_{\beta}: \mathcal{R}^{(n+1)}\to \mathcal{R}$ is the non-negative real valued function defined by
\begin{equation}\label{eq10}
\left|\left|\displaystyle\sum_{k=0}^{n} x_i\varepsilon ^i\right|\right|_{\beta}:=
\sqrt{\displaystyle\sum_{k=0}^{n} (n+1-i)\beta^i x_i^2\,}
\end{equation}
for $\displaystyle\sum_{k=0}^{n} x_i\varepsilon ^i\in \mathcal{R}^{(n+1)}$ with 
$x_i\in \mathcal{R}$ for $0\le i\le n$, then $||\,\,||_{\beta}$ makes  the truncated polynomial real algebra 
$\mathcal{R}^{(n+1)}$ into a normed algebra.
\end{proposition}

\medskip
\noindent
{\bf Proof} For convenience, we set $\alpha_i:=(n+1-i)\beta^i$ for $0\le i\le n$. Let $|\,\,|_{\tilde *}: \displaystyle\bigcup_{i=0}^{n}\mathcal{R}\varepsilon ^i\to \mathcal{R}$ be a map defined by 
\begin{equation}\label{eq11}
|x_i\varepsilon ^i|_{\tilde *}:=\sqrt{\alpha_i}\, |x_i|\quad\mbox{for $x_i\in \mathcal{R}$ and $0\le i\le n$}.
\end{equation}
For $x_i$ , $y_i$, $r\in \mathcal{R}$ and $0\le i\le n$, we clearly have
\begin{equation}\label{eq12}
\mbox{$|x_i\varepsilon ^i|_{\tilde *}\ge 0$, and $|x_i\varepsilon ^i|_{\tilde *}= 0$ if and only if $x_i\varepsilon ^i= 0$,}
\end{equation}
\begin{equation}\label{eq13}
\mbox{$|rx_i\varepsilon ^i|_{\tilde *}=|r|\,|x_i\varepsilon ^i|_{\tilde *}$}
\end{equation}
and
\begin{equation}\label{eq14}
\mbox{$|x_i\varepsilon ^i+y_i\varepsilon ^i|_{\tilde *}
=\sqrt{\alpha_i}\,|x_i+y_i|\le \sqrt{\alpha_i}\,(|x_i|+y_i|)\le 
|x_i\varepsilon ^i|_{\tilde *}+|y_i\varepsilon ^i|_{\tilde *}$}.
\end{equation}

We now prove that
\begin{equation}\label{eq15}
|x_i\varepsilon ^i \cdot y_j\varepsilon ^j|_{\tilde *}\le |x_i\varepsilon ^i|_{\tilde *} \cdot |y_j\varepsilon ^j|_{\tilde *}\quad\mbox{for $0\le i,\, j\le n$}.
\end{equation}

Since 
\begin{equation}\label{eq16}
|x_i\varepsilon ^i \cdot y_j\varepsilon ^j|_{\tilde *}=|(x_i y_j)\varepsilon ^{i+j}|_{\tilde *}
=\left\{\begin{array}{cl}\sqrt{\alpha_{i+j}}\,|x_i y_j|&
\mbox{if $i+j\le n$}\\0&\mbox{if $i+j> n$,}\end{array}\right.
\end{equation}
(\ref{eq15}) holds clearly if $i+j> n$. In the case where  $i+j\le n$,  we have
\begin{eqnarray*}
&&\alpha_{i}\alpha_{j}-\alpha_{i+j}=(n+1-i)\beta^i\,(n+1-j)\beta^j-\big(n+1-(i+j)\big)\beta^{i+j}\nonumber\\
&=&\beta^{i+j}\Big[(n+1)^2-(n+1)(i+j)+ij-\big(n+1-(i+j)\big)\Big]\nonumber\\
&=&\beta^{i+j}\Big[(n+1)\big(n+1-(i+j)\big)+ij-\big(n+1-(i+j)\big)\Big]\nonumber\\
&=&\beta^{i+j}\Big[n\big((n+1-(i+j)\big)+ij\Big]\ge 0\quad\mbox{for $i+j\le n$}
\end{eqnarray*}
or
\begin{equation}\label{eq17}
\alpha_{i+j}\le \alpha_{i}\alpha_{j}\quad\mbox{for $0\le i, , j\le n$ and $i+j\le n$.}
\end{equation}

It follows from (\ref{eq16}) and (\ref{eq17}) that
$$
|x_i\varepsilon ^i \cdot y_j\varepsilon ^j|_{\tilde *}\le \sqrt{\alpha_{i+j}}\,|x_i y_j|
\le \sqrt{\alpha_i}\,\sqrt{\alpha_j}\,|x_i|\,| y_j|
=\sqrt{\alpha_i}\,\,|x_i|\cdot\,\sqrt{\alpha_j}\,| y_j|,
$$
which proves that (\ref{eq15}) is also true if $i+j\le n$.

\medskip
By (\ref{eq12}), (\ref{eq13}), (\ref{eq14}) and (\ref{eq15}), $|\,\,|_{\tilde *}$ is a homogeneous norm on
 the truncated polynomial real algebra 
$\mathcal{R}^{(n+1)}$. By (\ref{eq10}) and (\ref{eq11}), we have
\begin{equation}\label{eq18}
\left|\left|\displaystyle\sum_{i=0}^{n} x_i\varepsilon ^i\right|\right|_{\beta}
=\left(\displaystyle\sum_{i=0}^{n}\big( |\alpha_i x_i|_{\tilde *}\big)^2\,\right)^{\frac12}
=\left|\left|\displaystyle\sum_{i=0}^{n} x_i\varepsilon ^i\right|\right|_{2,\,\tilde *}.
\end{equation}
It follows from (\ref{eq18}) and Proposition \ref{pr1.2} (iii) that  
$||\,\,||_{\beta}=||\,\,||_{2,\,\tilde *}$ is a norm on  the truncated polynomial real algebra
 $\mathcal{R}^{(n+1)}$.

\bigskip
In order to prove that  the truncated polynomial real  algebra $\mathcal{R}^{(n+1)}$ is a normed algebra
via  the norm $||\,\,||_{\beta}$, we need only to prove
\begin{equation}\label{eq19}
\left|\left|x\,y\right|\right|_{\beta}
\le \left|\left|x\right|\right|_{\beta}\,\left|\left|y\right|\right|_{\beta}\quad
\mbox{for $x$, $y\in \mathcal{R}^{(n+1)}$}.
\end{equation}

For $x=\displaystyle\sum_{i=0}^{n} x_i\varepsilon ^i$, we define 
\begin{equation}\label{eq20}
\phi(x)=
\left[\begin{array}{llllll}
x_0\beta^0&0&0&\cdots&0&0\\
x_1\beta&x_0\beta^0&0&\cdots&0&0\\
x_2\beta^2&x_1\beta&x_0\beta^0&\cdots&0&0\\
\vdots&\vdots&\vdots&\cdots&\vdots&\vdots\\
x_{n-1}\beta^{n-1}&x_{n-2}\beta^{n-2}&x_{n-3}\beta^{n-3}&\cdots&x_0\beta^0&0\\
x_n\beta^n&x_{n-1}\beta^{n-1}&x_{n-2}\beta^{n-2}&\cdots&x_1\beta&x_0\beta^0
\end{array}\right].
\end{equation}
Then  the map $\phi: \mathcal{R}^{(n+1)}\to {\bf M}_{n+1}(\mathcal{R})$ defined by
(\ref{eq20}) is  an injective algebra homomorphism. Using this algebra homomorphism $\phi$ and the matrix norm which makes  ${\bf M}_{n+1}(\mathcal{R})$ into a normed algebra, we get (\ref{eq19}).

\hfill\raisebox{1mm}{\framebox[2mm]{}}

\bigskip
\section {The  Definition of Automatic Differentiation }

Let $A$ be a unital associative real algebra.  For a non-empty set $S$,  we use 
$${\bf F}(S, A):= \{f\,|\,\mbox{$ f: S\to A$ is a function } \}$$
to denote the set of the functions from $S$ to $A$. For $f$, $g\in {\bf F}(S, A)$, $r\in \mathcal{R}$ and $x\in S$, we define
$$
(f+g)(x):=f(x)+g(x),\quad (rf)(x):=r\cdot f(x),\quad (f\cdot g)(x):=f(x)g(x).
$$
Then ${\bf F}(S, A)$ is a  unital associative real algebra with respect to the addition,  the scalar multiplication and the product above. The  identity $1_{{\bf F}(S, \, A)}$ of the algebra  ${\bf F}(S, A)$ is the constant function given by 
$$
1_{{\bf F}(S, \, A)} (x):=1\quad \mbox{for $x\in S$,}
$$
where $1$ is the identity of the  unital associative real algebra $A$.

For $x\in \mathcal{R}$, we define
$$
{\bf D}^n(x):=\left\{f\,\left|\begin{array}{l}
\mbox{$f$ is a real-valued function defined on an open interval}\\
\mbox{of real numbers and $f$ has the $n$-th derivative at $x$}
\end{array}\right.\right\}.
$$
Clearly, ${\bf D}^n(x)$ is a unital associative real algebra. 

\medskip
Using the notations above, we now  present our way of conceptualizing  automatic differentiation mathematically in the following definition. 

\medskip
\begin{definition}\label{def2.1} Let $A$ be a real-like algebra. A $3$-tuple  $\big(\Lambda,\, \Omega, \{\Gamma_{i}\}_{i=1}^n\big)$ 
consisting of  a map    
$\Lambda: \displaystyle\bigcup _{x\in \mathcal{R}}{\bf D}^n(x)\to  \displaystyle\bigcup_{S\subseteq A}{\bf F}(S, \,A)$, a map 
$\Omega: \mathcal{R}\to A$ and a family maps $\{\Gamma_{i}: A\to \mathcal{R}\,|\, 1\le i\le n\}$ is called the {\bf $n$-th order  automatic differentiation induced by $A$}  or  
{\bf  $A\big(\Lambda,\, \Omega, \{\Gamma_{i}\}_{i=1}^n\big)$-automatic differentiation} if the following four conditions are satisfied:
\begin{description} 
\item[(i)] For each $x\in \mathcal{R}$, there exists a subset $A_x\subseteq A$ such that $\Omega(x)\in A_x$,
$Im (\Lambda |{\bf D}^n(x))\subseteq {\bf F}(A_x, \,A)$ and the map  
$(\Lambda |{\bf D}^n(x)): {\bf D}^n(x))\to {\bf F}(A_x, \,A)$ is a real linear transformation;
\item[(ii)] $\Lambda$  preserves the product at $\Omega(x)$, which means 
\begin{equation}\label{eq32A}
\Lambda(f\cdot g)\big(\Omega(x)\big)=(\Lambda f)\big(\Omega(x)\big)\cdot (\Lambda g)\big(\Omega(x)\big)\quad
\mbox{for $f$, $g\in {\bf D}^n(x)$ };
\end{equation}
\item[(iii)]   $\Lambda$  preserves the composition  at $\Omega(x)$, which means 
\begin{equation}\label{eq32}
 \Lambda(f\circ g)\big(\Omega(x)\big)=\big(\Lambda(f)\circ \Lambda(g)\big)\big(\Omega(x)\big)\quad \mbox{for $f$, $g\in {\bf D}^n(x)$ };
\end{equation}
\item[(vi)]  The map $\Gamma_{i}: A\to \mathcal{R}$ for  $i\in \{1, 2, \dots , k\}$, which is called 
the {\bf $i$-th derivative map}, has the following property:
\begin{equation}\label{eq31}
\big(\Gamma_{i}\circ \Lambda(f)\big)\big(\Omega(x)\big)=f^{(i)}(x),
\end{equation}
where $f \in {\bf D}^n(x)$,  $x\in \mathcal{R}$ and $f^{(i)}$ is the $i$-th derivative of $f$.
\end{description}
\end{definition} 

\bigskip
Currently, the technique of evaluating high order derivatives in scientific computation community is the $n$-th order automatic differentiation induced by the  $(n+1)$-dimensional 
truncated polynomial real algebra $\mathcal{R}^{(n+1)}=\displaystyle\bigoplus_{i=0}^n \mathcal{R}\varepsilon ^i$, where the map  
$\Omega:=\Omega_{\theta_1}:  \mathcal{R}\to \mathcal{R}^{(n+1)}$ required in Definition \ref{def2.1} is given by
\begin{equation}\label{eq34A}
\Omega_{\theta_1}(x):=x+{\theta_1}\varepsilon\quad\mbox{for $x\in \mathcal{R}$.}
\end{equation}
The parameter ${\theta_1}$ in (\ref{eq34A}) is a non-zero real constant. Although the ways of defining the map 
$\Omega=\Omega_{\theta_1}$ required in Definition \ref{def2.1} are not unique,  the way given by  (\ref{eq34A}) is the only way which has appeared in the study of higher order automatic differentiation. In the remaining part of this paper, we will discuss how to replace the map $\Omega_{\theta_1}$ in (\ref{eq34A}) with the map 
$\Omega_{\theta_1, \theta_2, \dots , \theta_n}$ to introduce different ways of using automatic differentiation  technique to evaluate  higher order  derivatives, where $\theta_1\ne 0$, $\theta_2$, $\dots$ , $\theta_n$ are real constants 
and the map $\Omega_{\theta_1, \theta_2, \dots , \theta_n}: \mathcal{R}\to \mathcal{R}^{(n+1)}$ is given by
\begin{equation}\label{eq34B}
\Omega_{\theta_1, \theta_2, \dots , \theta_n}(x):=x+{\theta_1}\varepsilon+{\theta_2}\varepsilon^2+
\cdots +{\theta_n}\varepsilon^n\quad\mbox{for $x\in \mathcal{R}$.}
\end{equation}
For simplicity, we will use the $3$-rd order automatic differentiation induced by $\mathcal{R}^{(4)}$ as an example  to demonstrate the 
different ways of evaluating the first, the second and third derivatives of a differentiable function  exactly and simultaneously.
 
\medskip
Although the  $n$-th order automatic differentiation induced by the truncated polynomial real algebra $\mathcal{R}^{(n+1)}$
has been used in scientific computation for a long time,  its mathematical foundation has not been established  solidly.  For example, in order to prove the correctness of the  $n$-th order automatic differentiation induced by $\mathcal{R}^{(n+1)}$,  we must prove that the equation (\ref{eq32}) holds. To the best of our knowledge, the complete proof of  the equation (\ref{eq32}) for the  $n$-th order automatic differentiation  with $n\ge 2$ has not appeared in the study of automatic differentiation. The proof of  the equation (\ref{eq32}) is quite straightforward for  the  first order automatic differentiation as indicated in  Section 3.1.1 of \cite{GPRS}, but it takes some effort to prove  the equation (\ref{eq32})  for the  $n$-th order automatic differentiation  with $n\ge 2$.  The similar problems exists for the  mathematical foundation of multivariate high order automatic differentiation.  Let us mention some unsatisfying points about multivariate high order automatic differentiation.   In an attempt to develop multivariate high order automatic differentiation, we have to have the counterpart of the map  $\Lambda$ in  Definition \ref{def2.1}. Multivariate high order automatic differentiation is discussed in \cite{B} and \cite{H}. One way of defining  the counterpart of the map  $\Lambda$ in  Definition \ref{def2.1} in the context of   multivariate high order automatic differentiation is given by the map
$$
\hat{\hat{h}}_{n, N}=\hat{\hat{h}}: \,\Big(\mathcal{R}[X_1, \dots , X_n]/I_N\Big)^n_U\to \mathcal{R}[X_1, \dots , X_n]/I_N
$$
in (8.1) of \cite{H}, where $\mathcal{R}[X_1, \dots , X_n]/I_N$,  which is denoted by $_ND_n$ in \cite{B},  is the the truncated polynomial algebra and $I_N$ is the idea of the multi-variable polynomial ring  $\mathcal{R}[X_1, \dots , X_n]$ generated by the set 
$$
\left\{\left.X_1^{k_1}\cdots   X_n^{k_n}\,\right|\, \mbox{$k_1, \cdots , k_n$ are non-negative integers and $k_1+ \cdots +k_n>N$} \right\}.
$$
In order to  prove that the map $\hat{\hat{h}}_{n, N}=\hat{\hat{h}}$ works for evaluating high order partial derivatives of  multi-variable real-valued differentiable function using  automatic differentiation technique, it is really essential to prove 
\begin{equation}\label{eq33A}
\mbox{$\hat{\hat{h}}_{n, N}=\hat{\hat{h}}$ preserves the associative product}
\end{equation}
and
\begin{equation}\label{eq33B}
\mbox{$\hat{\hat{h}}_{n, N}=\hat{\hat{h}}$ preserves the composition of functions.}
\end{equation}
In \cite{H}, the proof of (\ref{eq33A}) has been omitted and the proof of (\ref{eq33B}) has not been discussed. In Section 2.3.1 of \cite{B}, although the necessity of proving (\ref{eq33B}) is mentioned, the discussion in Section 2.3.1 of \cite{B} is not a complete proof of 
\cite{B}.  The  proofs of  (\ref{eq33A}) and  (\ref{eq33B}) are not simple. From the work we have done, we see that the complete proofs of  (\ref{eq33A}) and  (\ref{eq33B})  mainly consist of three steps. First, we need to express the  map $\hat{\hat{h}}_{n, N}=\hat{\hat{h}}$  in a combinatorial form because the 
expression of  $\hat{\hat{h}}_{n, N}=\hat{\hat{h}}$ in  (8.1) of \cite{H} is neither suitable for proving   (\ref{eq33A}) nor  (\ref{eq33B}).
Next, we need to use the generalization of Leibniz Rule given in Proposition 6 of \cite{Ha} to prove   (\ref{eq33A}). Finally, we need to use the multivariate Faa Di Bruno formula introduced in \cite{CS}  to prove   (\ref{eq33B}). Hence, more work should be done to delete the 
unsatisfying points about multivariate high order automatic differentiation.

\bigskip
\section{The $3$-rd Order Automatic Differentiation}

We  now begin to explain how  to get  the $3$-rd order automatic differentiation induced by the   truncated polynomial real algebra $\mathcal{R}^{(4)}=\displaystyle\bigoplus_{i=0}^3\mathcal{R}\varepsilon^i$ . 

\medskip
For   $x\in \mathcal{R}$ and $f\in {\bf D}^3(x)$, we let
$$
A_x:=\{x+a_1\varepsilon +a_2\varepsilon^2+a_3\varepsilon^3\,|\, a_1, a_2, a_3\in \mathcal{R}\}
$$
and define the map
$\Lambda: {\bf D}^3(x)\to{\bf F}(A_x, \,\mathcal{R}^{(4)})$ by
$$
\Lambda (f)(x+a_1\varepsilon +a_2\varepsilon^2+a_3\varepsilon^3):=f(x)+a_1f'(x)\varepsilon+
$$
\begin{equation}\label{eq344}
+\Big(a_2f'(x)+\displaystyle\frac12 a_1^2f''(x)\Big)\varepsilon^2+\Big(a_3f'(x)+a_1a_2f''(x)+\displaystyle\frac16 a_1^3f'''(x)\Big)\varepsilon^3,
\end{equation}
where  $a_1$, $a_2$, $a_3\in \mathcal{R}$. 

Let $\alpha \ne 0$, $\beta$ and $\gamma$ be three real constants.  We define $\Omega_{\alpha, \beta, \gamma}: \mathcal{R}\to \mathcal{R}^{(4)}$  by
\begin{equation}\label{eq33}
\Omega_{\alpha, \beta, \gamma} (x):=x+\alpha \varepsilon +\beta\varepsilon^2+\gamma\varepsilon^3\quad\mbox{for $x\in \mathcal{R}$.}
\end{equation}

\medskip
The following theorem, which is the main theorem of this paper, presents the 
many ways  of doing automatic differentiation to compute the first, the second and the third  derivatives exactly and simultaneously.

\begin{proposition}\label{pr2.1}  Let $\alpha \ne 0$, $\beta$ and $\gamma$ be three real constants. If the map 
$\Lambda: \displaystyle\bigcup _{x\in \mathcal{R}}{\bf D}^n(x)\to  \displaystyle\bigcup_{S\subseteq \mathcal{R}^{(4)}}{\bf F}(S, \,\mathcal{R}^{(4)})$ 
and the map $\Omega_{\alpha, \beta, \gamma}: \mathcal{R}\to \mathcal{R}^{(4)}$  are defined by (\ref{eq344}) and (\ref{eq33}), then the $3$-tuple  $\big(\Lambda,\, \Omega_{\alpha, \beta, \gamma}, \{\Gamma_{i}\}_{i=1}^3\big)$ 
 is  the $3$-rd order automatic differentiation induced by the   truncated polynomial real algebra $\mathcal{R}^{(4)}$, where the 
$i$-th derivative map $\Gamma_{i}: \mathcal{R}^{(4)}\to \mathcal{R}$ for each $i\in \{1, 2, 3\}$ is defined by
\begin{equation}\label{eq37}
\left\{\begin{array}{l}
\Gamma_{1}(y):=\displaystyle\frac{1}{\alpha} y_1, \\\\
\Gamma_{2}(y):=\displaystyle\frac{2}{\alpha^2}y_2-\displaystyle\frac{2\beta}{\alpha^3}y_1,\\\\
 \Gamma_{3}(y):=\displaystyle\frac{6}{\alpha^3}y_3-\displaystyle\frac{12\beta}{\alpha^4}y_2+
\left(\displaystyle\frac{12\beta^2}{\alpha^5}-\displaystyle\frac{6\gamma}{\alpha^4}\right)y_1
\end{array}\right.
\end{equation}
for $y=y_0+y_1\varepsilon +y_2\varepsilon^2+y_3\varepsilon^3\in \mathcal{R}^{(4)}$ with
$y_0$, $y_1$, $y_2$, $y_3\in \mathcal{R}$.
\end{proposition}

\medskip
\noindent
{\bf Proof}  Throughout the proof of this proposition, we let  $x\in \mathcal{R}$ and $f$, $g\in {\bf D}^3(x)$.

\bigskip
First, by (\ref{eq344}), $(\Lambda |{\bf D}^3(x)): {\bf D}^3(x)\to{\bf F}(A_x, \,\mathcal{R}^{(4)})$  is clearly a real linear transformation. Hence,  the property (i)  Definition \ref{def2.1}  holds.

\bigskip
Next,   we have
\begin{equation}\label{eq47}
(fg)'=f'g+fg',\quad 
(fg)''=f''g+2f'g'+fg''
\end{equation}
and
\begin{equation}\label{eq48}
(fg)'''=f'''g+3f''g'+3f'g''+fg'''.
\end{equation}

Let $x+a_1\varepsilon +a_2\varepsilon^2+a_3\varepsilon^3$ be any element in  $\mathcal{R}^{(4)}$, where 
$x$, $a_1$, $a_2$ and $a_3\in \mathcal{R}$.
By (\ref{eq344}), (\ref{eq47}) and (\ref{eq48}), we have
\begin{eqnarray}\label{eq49}
&&\Lambda (fg)(x+a_1\varepsilon +a_2\varepsilon^2+a_3\varepsilon^3)=fg +a_1(fg)'\varepsilon +\nonumber\\
&&\quad +\Big[a_2(fg)'+\displaystyle\frac12 a_1^2(fg)''\Big]\varepsilon^2+\Big[a_3(fg)'+a_1a_2(fg)''+\displaystyle\frac16 a_1^3(fg)'''\Big]\varepsilon^3\nonumber\\
&=&fg +a_1(f'g+fg')\varepsilon+\Big[a_2(f'g+fg')+\displaystyle\frac12 a_1^2(f''g+2f'g'+fg''
)\Big]\varepsilon^2+\nonumber\\
&&\quad +\Big[a_3(f'g+fg')+a_1a_2(f''g+2f'g'+fg''
)+\nonumber\\
&&\quad +\displaystyle\frac16 a_1^3(f'''g+3f''g'+3f'g''+fg''')\Big]\varepsilon^3
\end{eqnarray}
and
\begin{eqnarray*}
&&\Big(\Lambda (f)\cdot \Lambda (g)\Big)(x+a_1\varepsilon +a_2\varepsilon^2+a_3\varepsilon^3)\nonumber\\
&=&\Lambda (f)(x+a_1\varepsilon +a_2\varepsilon^2+a_3\varepsilon^3)\cdot
\Lambda (g)(x+a_1\varepsilon +a_2\varepsilon^2+a_3\varepsilon^3)\nonumber\\
&=&\Big[f+a_1f'\varepsilon+\Big(a_2f'+\displaystyle\frac12 a_1^2f''\Big)\varepsilon^2+\Big(a_3f'+a_1a_2f''+\displaystyle\frac16 a_1^3f'''\Big)\varepsilon^3\Big]\cdot\nonumber\\
&&\quad
\cdot\Big[g+a_1g'\varepsilon+\Big(a_2g'+\displaystyle\frac12 a_1^2g''\Big)\varepsilon^2+\Big(a_3g'+a_1a_2g''+\displaystyle\frac16 a_1^3g'''\Big)\varepsilon^3\Big]\nonumber\\
&=&fg+(a_1f'\cdot g+f\cdot a_1g')\varepsilon+\nonumber\\
&&\quad +\Big[f\cdot \Big(\underbrace{a_2g'}_{1}+\underbrace{\displaystyle\frac12 a_1^2g''}_{2}\Big)+ \underbrace{a_1f'\cdot a_1g'}_{2}
+\Big(\underbrace{a_2f'}_{1}+\underbrace{\displaystyle\frac12 a_1^2f''}_{2}\Big)\cdot g\Big]\varepsilon^2+\nonumber\\
&&\quad +\Big[f\cdot \Big(\underbrace{a_3g'}_{3}+\underbrace{a_1a_2g''}_{4}+
\underbrace{\displaystyle\frac16 a_1^3g'''}_{5}\Big)
+a_1f'\cdot \Big(\underbrace{a_2g'}_{4}+\underbrace{\displaystyle\frac12 a_1^2g''}_{5}\Big)+\nonumber\\
\end{eqnarray*}
$$
+ \Big(\underbrace{a_2f'(x)}_{4}+\underbrace{\displaystyle\frac12 a_1^2f''}_{5}\Big)\cdot a_1g'+
\Big(\underbrace{a_3f'}_{3}+\underbrace{a_1a_2f''}_{4}+\underbrace{\displaystyle\frac16 a_1^3f'''}_{5}\Big)\cdot g\Big]\varepsilon^3
$$
\begin{eqnarray}\label{eq50}
&=&fg +a_1(f'g+fg')\varepsilon+\Big[a_2(f'g+fg')+\displaystyle\frac12 a_1^2(f''g+2f'g'+fg''
)\Big]\varepsilon^2+\nonumber\\
&&\quad +\Big[a_3(f'g+fg')+a_1a_2(f''g+2f'g'+fg''
)+\nonumber\\
&&\quad +\displaystyle\frac16 a_1^3(f'''g+3f''g'+3f'g''+fg''')\Big]\varepsilon^3.
\end{eqnarray}
Using (\ref{eq49}) and (\ref{eq50}), we get
\begin{equation}\label{eq51}
\Lambda (f\cdot g)=\Lambda (f)\cdot \Lambda (g)\quad\mbox{for $f$, $g\in  {\bf D}^3(x)$},
\end{equation}
which proves  that  $\Lambda$  preserves the product  in the algebra  ${\bf D}^3(x)$. In particular, the property (ii)  Definition \ref{def2.1}  holds.

\bigskip
Thirdly, it follows from (\ref{eq344}) and (\ref{eq33}) that
\begin{eqnarray*}
&&\big(\Lambda (f)\circ \Lambda (g)\big)\big(\Omega(x)\big)=\big(\Lambda (f)\circ \Lambda (g)\big)(x+\alpha \varepsilon +\beta\varepsilon^2+\gamma\varepsilon^3)\nonumber\\
&=&\Lambda (f)\big( \Lambda (g)(x+\alpha \varepsilon +\beta\varepsilon^2+\gamma\varepsilon^3)\big)
=\Lambda (f)\Big(g(x)+\alpha g'(x)\varepsilon+\nonumber\\
&&+\Big(\beta g'(x)+\displaystyle\frac12 \alpha^2g''(x)\Big)\varepsilon^2+
\Big(\gamma g'(x)+\alpha\beta g''(x)+\displaystyle\frac16 \alpha^3g'''(x)\Big)\varepsilon^3\Big)\nonumber\\
&=&f\big(g(x)\big)+\alpha g'(x)\cdot f'\big(g(x)\big)+\nonumber\\
&& +\Big[\Big(\beta g'(x)+\displaystyle\frac12 \alpha^2g''(x)\Big)\cdot f'\big(g(x)\big)+\displaystyle\frac12
\big(\alpha g'(x)\big)^2\cdot f''\big(g(x)\big)\Big]\varepsilon^2+\nonumber\\
&&+\Big[\Big(\gamma g'(x)+\alpha\beta g''(x)+\displaystyle\frac16 \alpha^3g'''(x)\Big)\cdot f'\big( g(x)\big)+\nonumber\\
\end{eqnarray*}
$$
+\alpha g'(x)\cdot \Big(\beta g'(x)+\displaystyle\frac12 \alpha^2g''(x)\Big)\cdot f''\big(g(x)\big)+
\displaystyle\frac16 \big(\alpha g'(x)\big)^3\cdot f'''\big(g(x)\big)\Big]\varepsilon^3+
$$
\begin{eqnarray*}\label{eq51A}
&=&(f\circ g)(x)+\alpha (f\circ g)'(x)\varepsilon+\Big(\beta  (f\circ g)'(x)+\displaystyle\frac12 \alpha^2 (f\circ g)''(x)\Big)\varepsilon^2+\\
&&+\Big(\gamma (f\circ g)'(x)+\alpha\beta (f\circ g)''(x)+\displaystyle\frac16 \alpha^3(f\circ g)'''(x)\Big)\varepsilon^3\\
&=&\big(\Lambda (f\circ g)\big)(x+\alpha \varepsilon +\beta\varepsilon^2+\gamma\varepsilon^3)
=\big(\Lambda (f\circ g)\big)\big(\Omega(x)\big),
\end{eqnarray*}
which  proves that  the property (iii)  Definition \ref{def2.1} holds.

\bigskip
Finally, using (\ref{eq34}) and (\ref{eq37}), we have
$$
\Big((\Lambda(f)\circ\Omega_{\alpha, \beta, \gamma}\Big)(x)=\Lambda(f)(x+\alpha\varepsilon +\beta\varepsilon^2+\gamma\varepsilon^3)
=f(x)+\underbrace{\alpha f'(x)}_{y_1}\varepsilon+
$$
\begin{equation}\label{eq34}
+\underbrace{\Big(\beta f'(x)+\displaystyle\frac12 \alpha^2f''(x)\Big)}_{y_2}\varepsilon^2+
\underbrace{\Big(\gamma f'(x)+\alpha\beta f''(x)+\displaystyle\frac16 \alpha^3f'''(x)\Big)}_{y_3}\varepsilon^3
\end{equation}
By (\ref{eq37}) and (\ref{eq50}), we get
$$
\Big(\Gamma_{1}\circ (\Lambda(f))\circ\Omega_{\alpha, \beta, \gamma}\Big)(x)=\displaystyle\frac{1}{\alpha}\cdot \alpha f'(x)=f'(x),
$$
$$
\Big(\Gamma_{2}\circ (\Lambda(f))\circ\Omega_{\alpha, \beta, \gamma}\Big)(x)=
\displaystyle\frac{2}{\alpha^2}\cdot \Big(\beta f'(x)+\displaystyle\frac12 \alpha^2f''(x)\Big)
-\displaystyle\frac{2\beta}{\alpha^3}\cdot \alpha f'(x)
=f''(x)
$$
and
$$
\Big(\Gamma_{3}\circ (\Lambda(f))\circ\Omega_{\alpha, \beta, \gamma}\Big)(x)=
\displaystyle\frac{6}{\alpha^3}\cdot \Big(\gamma f'(x)+\alpha\beta f''(x)+\displaystyle\frac16 \alpha^3f'''(x)\Big)+
$$
$$
-\displaystyle\frac{12\beta}{\alpha^4}\cdot \Big(\beta f'(x)+\displaystyle\frac12 \alpha^2f''(x)\Big)+
\left(\displaystyle\frac{12\beta^2}{\alpha^5}-\displaystyle\frac{6\gamma}{\alpha^4}\right)
\cdot \alpha f'(x)=f'''(x),
$$
which proves that   the property (vi)  Definition \ref{def2.1} holds..

\medskip
This completes the proof of Proposition \ref{pr2.1}.

\hfill\raisebox{1mm}{\framebox[2mm]{}}

\bigskip
For $f\in {\bf D}^3(c)$ with $c\in \mathcal{R}$, the function 
$\Lambda (f)\in {\bf F}(A_c, \,\mathcal{R}^{(4)})$, which is also denoted by $\overline{f}$, will be called the {\bf $\mathcal{R}^{(4)}$-extension} of $f$. The $\mathcal{R}^{(4)}$-extensions of some common elementary functions in ${\bf D}^3(c)$ are given as follows:
\begin{eqnarray*}
&\bullet&\displaystyle\frac{1}{x+a_1\varepsilon +a_2\varepsilon^2+a_3\varepsilon^3}=
\displaystyle\frac{1}{x}-\displaystyle\frac{a_1}{x^2}\varepsilon 
+\left(\displaystyle\frac{a_1^2}{x^3}-\displaystyle\frac{a_2}{x^2}\right)\varepsilon^2+\\
&&\quad +\left(-\displaystyle\frac{a_1^3}{x^4}+\displaystyle\frac{2a_1a_2}{x^3}
-\displaystyle\frac{a_3}{x^2}\right)\varepsilon^3\quad\mbox{for $0\ne x\in \mathcal{R}$}\\
&&\\
&\bullet&\overline{\exp} (x+a_1\varepsilon +a_2\varepsilon^2+a_3\varepsilon^3)=e^x +
a_1e^x \varepsilon +\left(\displaystyle\frac{a_1^2}{2}+a_2\right)e^x\varepsilon^2+\\
&&\quad +\left(\displaystyle\frac{a_1^3}{6}+a_1a_2+a_3\right)e^x\varepsilon^3\\
&&\\
&\bullet&\overline{\sin}\, (x+a_1\varepsilon +a_2\varepsilon^2+a_3\varepsilon^3)=\sin x +
(a_1\cos x)\varepsilon+\\
&&\quad +\left(-\displaystyle\frac{a_1^2}{2}\sin x+a_2\cos x\right)\varepsilon^2+
\left(-\displaystyle\frac{a_1^3}{6}\cos x-a_1a_2\sin x+a_3\cos x\right)\varepsilon^3\\
&&\\
&\bullet&\overline{\cos}\, (x+a_1\varepsilon +a_2\varepsilon^2+a_3\varepsilon^3)=\cos x -
(a_1\sin x)\varepsilon+\\
&&\quad +\left(-\displaystyle\frac{a_1^2}{2}\cos x-a_2\sin x\right)\varepsilon^2+
\left(\displaystyle\frac{a_1^3}{6}\sin x-a_1a_2\cos x-a_3\sin x\right)\varepsilon^3\\
&&\\
&\bullet&\overline{\ln}\,(x+a_1\varepsilon +a_2\varepsilon^2+a_3\varepsilon^3)=
\ln x+\displaystyle\frac{a_1}{x}\varepsilon 
+\left(-\displaystyle\frac{a_1^2}{2x^2}+\displaystyle\frac{a_2}{x}\right)\varepsilon^2+\\
&&\quad +\left(\displaystyle\frac{a_1^3}{3x^3}-\displaystyle\frac{a_1a_2}{x^2}
+\displaystyle\frac{a_3}{x}\right)\varepsilon^3\quad\mbox{for $0< x\in \mathcal{R}$}\\
&&\\
&\bullet&\overline{\arctan}\,(x+a_1\varepsilon +a_2\varepsilon^2+a_3\varepsilon^3)=
\arctan x+\displaystyle\frac{a_1}{1+x^2}\varepsilon +\\
&&\quad 
+\left(-\displaystyle\frac{a_1^2x}{(1+x^2)^2} +\displaystyle\frac{a_2}{1+x^2}\right)\varepsilon^2
 +\left(\displaystyle\frac{a_1^3(3x^2-1)}{3(1+x^2)^3}-\displaystyle\frac{2a_1a_2x}{(1+x^2)^2}
+\displaystyle\frac{a_3}{1+x^2}\right)\varepsilon^3
\end{eqnarray*}

\bigskip
Since every function $f\in {\bf D}^3(c)$ is formed in terms of some elementary functions in ${\bf D}^3(c)$,  we can evaluate the first derivative $f'(c)$, the second derivative $f''(c)$ and the third derivative $f'''(c)$ 
exactly and simultaneously by using   $\mathcal{R}^{(4)}\big(\Lambda,\, \Omega_{\alpha, \beta, \gamma}, \{\Gamma_{i}\}_{i=1}^3\big)$-automatic differentiation and the $\mathcal{R}^{(4)}$-extensions of those elementary functions which form the function $f\in {\bf D}^3(c)$. In the remaining part of this paper, 
for convenience,  $\mathcal{R}^{(4)}\big(\Lambda,\, \Omega_{\alpha, \beta, \gamma}, \{\Gamma_{i}\}_{i=1}^3\big)$-automatic differentiation introduce in Proposition \ref{pr1.2} is also denoted the {\bf   $\mathcal{R}^{(4)}_{\alpha, \beta, \gamma}$-automatic differentiation}. The different choices of the parameters $\alpha\ne 0$, $\beta$ and  $\gamma$ give different ways of doing automatic differentiation to get  the first, the second and the third derivatives of  a function in $ {\bf D}^3(c)$ exactly and simultaneously. 
We now do an example by choosing $(\alpha, \beta, \gamma)=(1, 0, 0)$ and $(1, 1, 1)$ to do the  $\mathcal{R}^{(4)}_{\alpha, \beta, \gamma}$-automatic differentiation in two ways.

\medskip
After denoting a  $\mathcal{R}^{(4)}$-number $x+a_1\varepsilon +a_2\varepsilon^2+a_3\varepsilon^3\in \mathcal{R}^{(4)}$  by a $4$-tuple $(x, a_1, a_2, a_3)$ of real numbers,  the algorithm, which compute the first derivative $f'(c)$, the second  derivative $f''(c)$ and the third derivative $f'''(c)$ of a function $f\in {\bf D}^3(c)$ at a real number $c$, can be written in a pseudocode as follows:

\bigskip
$\underline{\overline{\mbox{{\bf Algorithm}\quad $\mathcal{R}^{(4)}_{\alpha, \beta, \gamma}$-automatic differentiation}}}$

\begin{itemize}
\item {\bf Input:} A real number $c$ and a function $f\in {\bf D}^3(c)$.
\item {\bf Output:} A $3$-tuple $(z_1, z_2, z_3)$ of real numbers, where $z_1=f'(c)$, $z_2=f''(c)$ and $z_3=f'''(c)$.
\end{itemize}

\begin{enumerate}
\item Start
\item Get the $\mathcal{R}^{(4)}$-extension $\overline{f}: \mathcal{R}^{(4)}\to \mathcal{R}^{(4)}$.
\item Compute the value of the function $\overline{f}$ at $(c, \alpha, \beta, \gamma)$ to get the $4$-tuple 
$(f(c), y_1, y_2, y_3)$ of real numbers.
\item Compute  $z_1=\displaystyle\frac{1}{\alpha} y_1$.
\item Compute  $z_2=\displaystyle\frac{2}{\alpha^2}y_2-\displaystyle\frac{2\beta}{\alpha^3}y_1$.
\item Compute  $z_3=\displaystyle\frac{6}{\alpha^3}y_3-\displaystyle\frac{12\beta}{\alpha^4}y_2+
\left(\displaystyle\frac{12\beta^2}{\alpha^5}-\displaystyle\frac{6\gamma}{\alpha^4}\right)y_1$.
\item Display  $(z_1, z_2, z_3)$
\item Stop
\end{enumerate}

\bigskip
As an example, let $f$ be a function defined by
\begin{equation}\label{eq38}
f(x)=(\ln x)\,\cos\Big(\displaystyle\frac{1}{x^2}\Big),
\end{equation}
which is formed in terms of the elementary functions $x^2$, $\displaystyle\frac{1}{x}$, $\ln x$ and  $\cos x$.
We now use both  $\mathcal{R}^{(4)}_{1, 0, 0}$-automatic differentiation and  $\mathcal{R}^{(4)}_{1, 1, 1}$-automatic differentiation to compute the first, the second and the third derivatives of the function $f$ defined by (\ref{eq38})  at $x\ne 0$. 
One thing we would like to emphasize is that although the symbols appear everywhere, the technique of evaluating derivatives in our examples are not symbolic differentiation technique.

\bigskip
{\bf $\bullet$ The example of doing  $\mathcal{R}^{(4)}_{1, 0, 0}$-automatic differentiation}. 

Using (\ref{eq33}) and the $\mathcal{R}^{(4)}$-extensions of the elementary functions $x^2$, $\displaystyle\frac{1}{x}$, $\ln x$ and  $\cos x$, we have
\begin{eqnarray*}
&&(\Lambda(f))\circ\Omega_{1, 0, 0}\Big)(x)=\overline{f}\Big(\Omega_{1, 0, 0} (x)\Big)=\overline{f}(x+\varepsilon )\\
&=&\overline{\ln} (x+\varepsilon )\cdot
\overline{\cos} \left(\displaystyle\frac{1}{(x+\varepsilon )^2}\right)
=\overline{\ln} (x+\varepsilon )\cdot
\overline{\cos} \left(\displaystyle\frac{1}{x^2+2x\varepsilon +\varepsilon^2}\right)\\
&=&\overline{\ln} (x+\varepsilon )\cdot
\overline{\cos} \left[\displaystyle\frac{1}{x^2}-\displaystyle\frac{2x}{(x^2)^2}\varepsilon \right.+\\
&&\quad +\left.\left(\displaystyle\frac{(2x)^2}{(x^2)^3}-\displaystyle\frac{1}{(x^2)^2}\right)\varepsilon^2+
\left(-\displaystyle\frac{(2x)^3}{(x^2)^4}+\displaystyle\frac{2\cdot 2x\cdot 1}{(x^2)^3}\right)\varepsilon^3\right]\\
&=&\overline{\ln} (x+\varepsilon )\cdot
 \overline{\cos} \left(\displaystyle\frac{1}{x^2}-\displaystyle\frac{2}{x^3}\varepsilon +\displaystyle\frac{3}{x^4}\varepsilon^2
-\displaystyle\frac{4}{x^5}\varepsilon^3\right)\\
&=&\left(\ln x+\displaystyle\frac{1}{x}\varepsilon+
-\displaystyle\frac{1}{2x^2}\varepsilon^2+
\displaystyle\frac{1}{3x^3}\varepsilon^3\right)\cdot \left\{\cos\left(\displaystyle\frac{1}{x^2}\right)-\left(-\displaystyle\frac{2}{x^3}\right)
\sin \left(\displaystyle\frac{1}{x^2}\right)\varepsilon+\right.\\
&&\quad + \left[-\displaystyle\frac{1}{2}\left(\displaystyle\frac{-2}{x^3}\right)^2\cos\left(\displaystyle\frac{1}{x^2}\right)-\displaystyle\frac{3}{x^4}\sin\left(\displaystyle\frac{1}{x^2}\right)\right]\varepsilon^2+
\end{eqnarray*}
$$
 +\left.\left[\displaystyle\frac{1}{6} \left(\displaystyle\frac{-2}{x^2}\right)^3
\sin\left(\displaystyle\frac{1}{x^2}\right)-\left(-\displaystyle\frac{2}{x^3}\right)
\displaystyle\frac{3}{x^4}\cos\left(\displaystyle\frac{1}{x^2}\right)-
\left(-\displaystyle\frac{4}{x^5}\right)\sin\left(\displaystyle\frac{1}{x^2}\right)
\right] \varepsilon^3\right\}
$$
\begin{eqnarray*}
&=&\left(\ln x+\displaystyle\frac{1}{x}\varepsilon+
-\displaystyle\frac{1}{2x^2}\varepsilon^2+
\displaystyle\frac{1}{3x^3}\varepsilon^3\right)\cdot \left\{\cos\left(\displaystyle\frac{1}{x^2}\right)+\displaystyle\frac{2}{x^3}
\sin \left(\displaystyle\frac{1}{x^2}\right)\varepsilon+\right.\\
&&\quad + \left[-\displaystyle\frac{2}{x^6}\cos\left(\displaystyle\frac{1}{x^2}\right)-\displaystyle\frac{3}{x^4}\sin\left(\displaystyle\frac{1}{x^2}\right)\right]\varepsilon^2+\\
&& +\left.\left[-\displaystyle\frac{4}{3x^9}
\sin\left(\displaystyle\frac{1}{x^2}\right)+\displaystyle\frac{6}{x^7}
\cos\left(\displaystyle\frac{1}{x^2}\right)+
\displaystyle\frac{4}{x^5}\sin\left(\displaystyle\frac{1}{x^2}\right)
\right] \varepsilon^3\right\}\\
&=&(\ln x)\,\cos\Big(\displaystyle\frac{1}{x^2}\Big)+\widetilde{y_1}\,\varepsilon +\widetilde{y_2}\,\varepsilon ^2+
\widetilde{y_3}\,\varepsilon ^3,
\end{eqnarray*}
where
\begin{equation}\label{eq39'}
\widetilde{y_1}=(\ln x)\cdot \displaystyle\frac{2}{x^3}\sin \Big(\displaystyle\frac{1}{x^2}\Big)
+\displaystyle\frac{1}{x}\cdot\cos \Big(\displaystyle\frac{1}{x^2}\Big),
\end{equation}
\begin{eqnarray}\label{eq40'}
&&\widetilde{y_2}=(\ln x)\cdot \left[-\displaystyle\frac{2}{x^6}\cos\left(\displaystyle\frac{1}{x^2}\right)-\displaystyle\frac{3}{x^4}\sin\left(\displaystyle\frac{1}{x^2}\right)\right]
+\displaystyle\frac{1}{x}\cdot \displaystyle\frac{2}{x^3}
\sin \left(\displaystyle\frac{1}{x^2}\right)+\nonumber\\
&&\quad\quad\quad 
 -\displaystyle\frac{1}{2x^2}\cdot \cos\left(\displaystyle\frac{1}{x^2}\right)\nonumber\\
&=&-\displaystyle\frac{2}{x^6}(\ln x)\cos\displaystyle\frac{1}{x^2}
-\displaystyle\frac{3}{x^4}(\ln x)\sin\displaystyle\frac{1}{x^2}
+\displaystyle\frac{2}{x^4}
\sin \displaystyle\frac{1}{x^2}
 -\displaystyle\frac{1}{2x^2}\cos\displaystyle\frac{1}{x^2}
\end{eqnarray}
and
$$
\widetilde{y_3}=(\ln x)\cdot \left[-\displaystyle\frac{4}{3x^9}
\sin\left(\displaystyle\frac{1}{x^2}\right)+\displaystyle\frac{6}{x^7}
\cos\left(\displaystyle\frac{1}{x^2}\right)+
\displaystyle\frac{4}{x^5}\sin\left(\displaystyle\frac{1}{x^2}\right)
\right]+
$$
$$
+\displaystyle\frac{1}{x}\cdot \left[-\displaystyle\frac{2}{x^6}\cos\left(\displaystyle\frac{1}{x^2}\right)-\underbrace{\displaystyle\frac{3}{x^4}\sin\left(\displaystyle\frac{1}{x^2}\right)}_{1}\right]
-\underbrace{\displaystyle\frac{1}{2x^2}\cdot \displaystyle\frac{2}{x^3}
\sin \left(\displaystyle\frac{1}{x^2}\right)}_{1}+\displaystyle\frac{1}{3x^3}\cdot \cos\left(\displaystyle\frac{1}{x^2}\right)
$$
\begin{eqnarray}\label{eq41'}
&=&-\displaystyle\frac{4}{3x^9}(\ln x)
\sin\left(\displaystyle\frac{1}{x^2}\right)+\displaystyle\frac{6}{x^7}(\ln x)
\cos\left(\displaystyle\frac{1}{x^2}\right)+
\displaystyle\frac{4}{x^5}(\ln x)\sin\left(\displaystyle\frac{1}{x^2}\right)+\nonumber\\
&&\qquad -\displaystyle\frac{2}{x^7}\cos\left(\displaystyle\frac{1}{x^2}\right)
-\displaystyle\frac{4}{x^5}\sin\left(\displaystyle\frac{1}{x^2}\right)
+\displaystyle\frac{1}{3x^3} \cos\left(\displaystyle\frac{1}{x^2}\right).
\end{eqnarray}

Using the $i$-th derivative map $\Gamma_{i}$ defined by (\ref{eq37}) for $i\in \{1, 2, 3\}$ and  $(\alpha, \beta, \gamma)=(1, 0, 0)$, we get from (\ref{eq39'}), (\ref{eq40'}) and (\ref{eq41'}) that
\begin{equation}\label{eq42'}
\Gamma_{1}(y)=\widetilde{y_1}=\displaystyle\frac{2}{x^3}(\ln x)\sin \Big(\displaystyle\frac{1}{x^2}\Big)
+\displaystyle\frac{1}{x}\cos \Big(\displaystyle\frac{1}{x^2}\Big)=f'(x),
\end{equation}
\begin{eqnarray}\label{eq43'}
&&\Gamma_{2}(y)=2\widetilde{y_2}=-\displaystyle\frac{4}{x^6}(\ln x)\cos\displaystyle\frac{1}{x^2}
-\displaystyle\frac{6}{x^4}(\ln x)\sin\displaystyle\frac{1}{x^2}
+\displaystyle\frac{4}{x^4}
\sin \displaystyle\frac{1}{x^2}+\nonumber\\
&&\qquad \qquad
 -\displaystyle\frac{1}{x^2}\cos\displaystyle\frac{1}{x^2}=f''(x)
\end{eqnarray}
and
$$
\Gamma_{3}(y)=6\widetilde{y_3}=-\displaystyle\frac{8}{x^9}(\ln x)
\sin\left(\displaystyle\frac{1}{x^2}\right)+\displaystyle\frac{36}{x^7}(\ln x)
\cos\left(\displaystyle\frac{1}{x^2}\right)+\displaystyle\frac{24}{x^5}(\ln x)\sin\left(\displaystyle\frac{1}{x^2}\right)+
$$
\begin{equation}\label{eq44'}
 -\displaystyle\frac{12}{x^7}\cos\left(\displaystyle\frac{1}{x^2}\right)
-\displaystyle\frac{24}{x^5}\sin\left(\displaystyle\frac{1}{x^2}\right)
+\displaystyle\frac{2}{x^3} \cos\left(\displaystyle\frac{1}{x^2}\right)
=f'''(x).
\end{equation}

\medskip
By (\ref{eq42'}), (\ref{eq43'}) and (\ref{eq44'}), we get the first derivative $f'(x)$, the second derivative $f''(x)$ and the third derivative $f'''(x)$ of 
$f(x)=(\ln x)\,\cos\Big(\displaystyle\frac{1}{x^2}\Big)$ by doing 
$\mathcal{R}^{(4)}_{1, 0, 0}$-automatic differentiation and without using any differentiation rule in calculus. 

\bigskip
{\bf $\bullet$ The example of doing $\mathcal{R}^{(4)}_{1, 1, 1}$-automatic differentiation}. 

Using (\ref{eq33}) and the $\mathcal{R}^{(4)}$-extensions of the elementary functions $x^2$, $\displaystyle\frac{1}{x}$, $\ln x$ and  $\cos x$, we have
\begin{eqnarray*}
&&(\Lambda(f))\circ\Omega_{1, 1, 1}\Big)(x)=\overline{f}\Big(\Omega_{1, 1, 1} (x)\Big)=\overline{f}(x+\varepsilon +\varepsilon^2+\varepsilon^3)\\
&=&\overline{\ln} (x+\varepsilon +\varepsilon^2+\varepsilon^3)\cdot
\overline{\cos} \left(\displaystyle\frac{1}{(x+\varepsilon +\varepsilon^2+\varepsilon^3)^2}\right)\\
&=&\ln (x+\varepsilon +\varepsilon^2+\varepsilon^3)\cdot
\cos \left(\displaystyle\frac{1}{x^2+2x\varepsilon +(2x+1)\varepsilon^2+(2x+2)\varepsilon^3}\right)\\
&=&\left[\ln x+\displaystyle\frac{1}{x}\varepsilon+
\left(-\displaystyle\frac{1}{2x^2}+\displaystyle\frac{1}{x}\right)\varepsilon^2+
\left(\displaystyle\frac{1}{3x^3}-\displaystyle\frac{1}{x^2}+\displaystyle\frac{1}{x}
\right)\varepsilon^3\right]\cdot \cos\left[\displaystyle\frac{1}{x^2}
-\displaystyle\frac{2x}{(x^2)^2}\varepsilon\right.\\
&&\quad \left.+\left(\displaystyle\frac{(2x)^2}{(x^2)^3}-\displaystyle\frac{2x+1}{(x^2)^2}\right)\varepsilon^2+
\left(-\displaystyle\frac{(2x)^3}{(x^2)^4}+\displaystyle\frac{2(2x)(2x+1)}{(x^2)^3}-
\displaystyle\frac{2x+2}{(x^2)^2}\right)\varepsilon^3\right]\\
&=&\left[\ln x+\displaystyle\frac{1}{x}\varepsilon+
\left(-\displaystyle\frac{1}{2x^2}+\displaystyle\frac{1}{x}\right)\varepsilon^2+
\left(\displaystyle\frac{1}{3x^3}-\displaystyle\frac{1}{x^2}+\displaystyle\frac{1}{x}
\right)\varepsilon^3\right]\cdot \cos\left[\displaystyle\frac{1}{x^2}
-\displaystyle\frac{2}{x^3}\varepsilon+\right.\\
&&\quad \left.+\left(\displaystyle\frac{4}{x^4}-\displaystyle\frac{2}{x^3}
-\displaystyle\frac{1}{x^4}
\right)\varepsilon^2+
\left(-\displaystyle\frac{4}{x^5}+\displaystyle\frac{6}{x^4}-
\displaystyle\frac{2}{x^3}\right)\varepsilon^3\right]\\
&=&\left[\ln x+\displaystyle\frac{1}{x}\varepsilon+
\left(-\displaystyle\frac{1}{2x^2}+\displaystyle\frac{1}{x}\right)\varepsilon^2+
\left(\displaystyle\frac{1}{3x^3}-\displaystyle\frac{1}{x^2}+\displaystyle\frac{1}{x}
\right)\varepsilon^3\right]\cdot \left\{
\cos \Big(\displaystyle\frac{1}{x^2}\Big)+\right.\\
&&\quad -\Big(-\displaystyle\frac{2}{x^3}\Big)\sin \Big(\displaystyle\frac{1}{x^2}\Big)\varepsilon+
\left[-\displaystyle\frac{1}{2}\Big(-\displaystyle\frac{2}{x^3}\Big)^2
\cos \Big(\displaystyle\frac{1}{x^2}\Big)-\Big(\displaystyle\frac{3}{x^4}\right.+\\
&&\quad -\left.\displaystyle\frac{2}{x^3}\Big)\sin \Big(\displaystyle\frac{1}{x^2}\Big)\right]\varepsilon^2+
\left[\displaystyle\frac{1}{6}\Big(-\displaystyle\frac{2}{x^3}\Big)^3
\sin \Big(\displaystyle\frac{1}{x^2}\Big)\right.+\\
&&\quad -\left.\left.\Big(-\displaystyle\frac{2}{x^3}\Big)
\Big(\displaystyle\frac{3}{x^4}-\displaystyle\frac{2}{x^3}\Big)
\cos \Big(\displaystyle\frac{1}{x^2}\Big)
-\Big(-\displaystyle\frac{4}{x^5}+\displaystyle\frac{6}{x^4}
-\displaystyle\frac{2}{x^3}\Big)\sin \Big(\displaystyle\frac{1}{x^2}\Big)\right]\varepsilon^3\right\}\\
&=&\left[\ln x+\displaystyle\frac{1}{x}\varepsilon+
\left(-\displaystyle\frac{1}{2x^2}+\displaystyle\frac{1}{x}\right)\varepsilon^2+
\left(\displaystyle\frac{1}{3x^3}-\displaystyle\frac{1}{x^2}+\displaystyle\frac{1}{x}
\right)\varepsilon^3\right]\cdot \left\{
\cos \Big(\displaystyle\frac{1}{x^2}\Big)+\right.\\
&&\quad +\displaystyle\frac{2}{x^3}
\sin \Big(\displaystyle\frac{1}{x^2}\Big)\varepsilon+
+\left[-\displaystyle\frac{2}{x^6}
\cos \Big(\displaystyle\frac{1}{x^2}\Big)+\Big(\displaystyle\frac{2}{x^3} -\displaystyle\frac{3}{x^4}\Big)\sin \Big(\displaystyle\frac{1}{x^2}\Big)\right]\varepsilon^2+
\end{eqnarray*}
$$+
\left.\left[-\displaystyle\frac{4}{3x^9}
\sin \Big(\displaystyle\frac{1}{x^2}\Big)+
\Big(\displaystyle\frac{6}{x^7}-\displaystyle\frac{4}{x^6}\Big)
\cos \Big(\displaystyle\frac{1}{x^2}\Big)
+\Big(\displaystyle\frac{4}{x^5}-\displaystyle\frac{6}{x^4}
+\displaystyle\frac{2}{x^3}\Big)\sin \Big(\displaystyle\frac{1}{x^2}\Big)\right]\varepsilon^3\right\}
$$
$$
=(\ln x)\,\cos\Big(\displaystyle\frac{1}{x^2}\Big)+y_1\varepsilon +y_2\varepsilon ^2+
y_3\varepsilon ^3,\qquad\qquad\qquad\qquad\qquad\qquad\qquad
$$
where
\begin{equation}\label{eq39}
y_1=\displaystyle\frac{2}{x^3}(\ln x)\sin \Big(\displaystyle\frac{1}{x^2}\Big)
+\displaystyle\frac{1}{x}\cos \Big(\displaystyle\frac{1}{x^2}\Big),
\end{equation}
\begin{eqnarray}\label{eq40}
y_2&=&(\ln x)\cdot\left[-\displaystyle\frac{2}{x^6}
\cos \Big(\displaystyle\frac{1}{x^2}\Big)+\Big(\displaystyle\frac{2}{x^3} -\displaystyle\frac{3}{x^4}\Big)\sin \Big(\displaystyle\frac{1}{x^2}\Big)\right]+\nonumber\\
&&\quad +\displaystyle\frac{1}{x}\cdot \displaystyle\frac{2}{x^3}
\sin \Big(\displaystyle\frac{1}{x^2}\Big)+
\left(-\displaystyle\frac{1}{2x^2}+\displaystyle\frac{1}{x}\right)\cdot 
\cos \Big(\displaystyle\frac{1}{x^2}\Big)\nonumber\\
&=&-\displaystyle\frac{2}{x^6}(\ln x)
\cos \Big(\displaystyle\frac{1}{x^2}\Big)+\Big(\displaystyle\frac{2}{x^3} -\displaystyle\frac{3}{x^4}\Big)(\ln x)
\sin \Big(\displaystyle\frac{1}{x^2}\Big)+\nonumber\\
&&\quad +\displaystyle\frac{2}{x^4}\sin \Big(\displaystyle\frac{1}{x^2}\Big)+
\left(-\displaystyle\frac{1}{2x^2}+\displaystyle\frac{1}{x}\right)
\cos \Big(\displaystyle\frac{1}{x^2}\Big)
\end{eqnarray}
and
$$
y_3=(\ln x)\cdot \left[-\displaystyle\frac{4}{3x^9}
\sin \Big(\displaystyle\frac{1}{x^2}\Big)+
\Big(\displaystyle\frac{6}{x^7}-\displaystyle\frac{4}{x^6}\Big)
\cos \Big(\displaystyle\frac{1}{x^2}\Big)\right.+
$$
$$
+\left.\Big(\displaystyle\frac{4}{x^5}-\displaystyle\frac{6}{x^4}
+\displaystyle\frac{2}{x^3}\Big)\sin \Big(\displaystyle\frac{1}{x^2}\Big)\right]+
\displaystyle\frac{1}{x}\cdot \left[-\displaystyle\frac{2}{x^6}
\cos \Big(\displaystyle\frac{1}{x^2}\Big)+\Big(\displaystyle\frac{2}{x^3} -\displaystyle\frac{3}{x^4}\Big)\sin \Big(\displaystyle\frac{1}{x^2}\Big)\right]+
$$
$$
\left(-\displaystyle\frac{1}{2x^2}+\displaystyle\frac{1}{x}\right)\cdot 
\displaystyle\frac{2}{x^3}
\sin \Big(\displaystyle\frac{1}{x^2}\Big)+
\left(\displaystyle\frac{1}{3x^3}-\displaystyle\frac{1}{x^2}+\displaystyle\frac{1}{x}
\right)\cdot \cos \Big(\displaystyle\frac{1}{x^2}\Big)
$$
or
$$
y_3=-\displaystyle\frac{4}{3x^9}(\ln x)\sin \Big(\displaystyle\frac{1}{x^2}\Big)+
\Big(\displaystyle\frac{6}{x^7}-\displaystyle\frac{4}{x^6}\Big)(\ln x)
\cos \Big(\displaystyle\frac{1}{x^2}\Big)+
$$
$$
+\Big(\displaystyle\frac{4}{x^5}-\displaystyle\frac{6}{x^4}
+\displaystyle\frac{2}{x^3}\Big)(\ln x)\sin \Big(\displaystyle\frac{1}{x^2}\Big)
-\displaystyle\frac{2}{x^7}
\cos \Big(\displaystyle\frac{1}{x^2}\Big)+
\Big(\displaystyle\frac{2}{x^4} -\displaystyle\frac{3}{x^5}\Big)\sin \Big(\displaystyle\frac{1}{x^2}\Big)+
$$
\begin{equation}\label{eq41}
+\left(-\displaystyle\frac{1}{x^5}+\displaystyle\frac{2}{x^4}\right)\sin \Big(\displaystyle\frac{1}{x^2}\Big)+
\left(\displaystyle\frac{1}{3x^3}-\displaystyle\frac{1}{x^2}+\displaystyle\frac{1}{x}
\right)\cos \Big(\displaystyle\frac{1}{x^2}\Big).
\end{equation}

Using the $i$-th derivative map $\Gamma_{i}$ defined by (\ref{eq37}) for $i\in \{1, 2, 3\}$ and  $(\alpha, \beta, \gamma)=(1, 1, 1)$, we get from (\ref{eq39}), (\ref{eq40}) and (\ref{eq41}) that
\begin{equation}\label{eq42}
\Gamma_{1}(y)=y_1=\displaystyle\frac{2}{x^3}(\ln x)\sin \Big(\displaystyle\frac{1}{x^2}\Big)
+\displaystyle\frac{1}{x}\cos \Big(\displaystyle\frac{1}{x^2}\Big)=f'(x),
\end{equation}
\begin{eqnarray}\label{eq43}
&&\Gamma_{2}(y)=2y_2-2y_1\nonumber\\
&=&
2\left[-\displaystyle\frac{2}{x^6}(\ln x)
\cos \Big(\displaystyle\frac{1}{x^2}\Big)+\Big(\displaystyle\frac{2}{x^3} -\displaystyle\frac{3}{x^4}\Big)(\ln x)
\sin \Big(\displaystyle\frac{1}{x^2}\Big)\right.+\nonumber\\
&&\quad +\left.\displaystyle\frac{2}{x^4}\sin \Big(\displaystyle\frac{1}{x^2}\Big)+
\left(-\displaystyle\frac{1}{2x^2}+\displaystyle\frac{1}{x}\right)
\cos \Big(\displaystyle\frac{1}{x^2}\Big)
\right]+\nonumber\\
&&\quad -2\left[\displaystyle\frac{2}{x^3}(\ln x)\sin \Big(\displaystyle\frac{1}{x^2}\Big)
+\displaystyle\frac{1}{x}\cos \Big(\displaystyle\frac{1}{x^2}\Big)\right]\nonumber\\
&=&-\displaystyle\frac{4}{x^6}(\ln x)
\cos \Big(\displaystyle\frac{1}{x^2}\Big)+\Big(\underbrace{\displaystyle\frac{4}{x^3}}_{1} -\displaystyle\frac{6}{x^4}\Big)(\ln x)
\sin \Big(\displaystyle\frac{1}{x^2}\Big) +\displaystyle\frac{4}{x^4}\sin \Big(\displaystyle\frac{1}{x^2}\Big)+\nonumber\\
&&\quad +
\left(-\displaystyle\frac{1}{x^2}+\underbrace{\displaystyle\frac{2}{x}}_{2}\right)
\cos \Big(\displaystyle\frac{1}{x^2}\Big)
-\underbrace{\displaystyle\frac{4}{x^3}(\ln x)\sin \Big(\displaystyle\frac{1}{x^2}\Big)}_{1}
-\underbrace{\displaystyle\frac{2}{x}\cos \Big(\displaystyle\frac{1}{x^2}}_{2}\Big)\nonumber\\
&=&-\displaystyle\frac{4}{x^6}(\ln x)
\cos \Big(\displaystyle\frac{1}{x^2}\Big) -\displaystyle\frac{1}{x^2}\cos \Big(\displaystyle\frac{1}{x^2}\Big)-\displaystyle\frac{6}{x^4}(\ln x)
\sin \Big(\displaystyle\frac{1}{x^2}\Big) +\nonumber\\
&&\quad +\displaystyle\frac{4}{x^4}\sin \Big(\displaystyle\frac{1}{x^2}\Big)=f''(x)
\end{eqnarray}
and
\begin{eqnarray*}
&&\Gamma_{3}(y)=6y_3-12y_2+6y_1\nonumber\\
&=&6\left[-\displaystyle\frac{4}{3x^9}(\ln x)\sin \Big(\displaystyle\frac{1}{x^2}\Big)+
\Big(\displaystyle\frac{6}{x^7}-\displaystyle\frac{4}{x^6}\Big)(\ln x)
\cos \Big(\displaystyle\frac{1}{x^2}\Big)\right.+\nonumber\\
\end{eqnarray*}
$$
+\Big(\displaystyle\frac{4}{x^5}-\displaystyle\frac{6}{x^4}
+\displaystyle\frac{2}{x^3}\Big)(\ln x)\sin \Big(\displaystyle\frac{1}{x^2}\Big)
-\displaystyle\frac{2}{x^7}
\cos \Big(\displaystyle\frac{1}{x^2}\Big)+\Big(\displaystyle\frac{2}{x^4} -\displaystyle\frac{3}{x^5}\Big)\sin \Big(\displaystyle\frac{1}{x^2}\Big)+
$$
\begin{eqnarray*}
&&\quad +\left.\left(-\displaystyle\frac{1}{x^5}+\displaystyle\frac{2}{x^4}\right)\sin \Big(\displaystyle\frac{1}{x^2}\Big)+
\left(\displaystyle\frac{1}{3x^3}-\displaystyle\frac{1}{x^2}+\displaystyle\frac{1}{x}
\right)\cos \Big(\displaystyle\frac{1}{x^2}\Big)\right]+\nonumber\\
&&\quad -12\left[-\displaystyle\frac{2}{x^6}(\ln x)
\cos \Big(\displaystyle\frac{1}{x^2}\Big)+\Big(\displaystyle\frac{2}{x^3} -\displaystyle\frac{3}{x^4}\Big)(\ln x)
\sin \Big(\displaystyle\frac{1}{x^2}\Big)\right.+\nonumber
\end{eqnarray*}
$$ 
+\left.\displaystyle\frac{2}{x^4}\sin \Big(\displaystyle\frac{1}{x^2}\Big)+
\left(-\displaystyle\frac{1}{2x^2}+\displaystyle\frac{1}{x}\right)
\cos \Big(\displaystyle\frac{1}{x^2}\Big)\right] +6\left[\displaystyle\frac{2}{x^3}(\ln x)\sin \Big(\displaystyle\frac{1}{x^2}\Big)
+\displaystyle\frac{1}{x}\cos \Big(\displaystyle\frac{1}{x^2}\Big)\right]
$$
\begin{eqnarray*}
&=&-\displaystyle\frac{8}{x^9}(\ln x)\sin \Big(\displaystyle\frac{1}{x^2}\Big)+
\Big(\displaystyle\frac{36}{x^7}-\underbrace{\displaystyle\frac{24}{x^6}}_{1}\Big)(\ln x)
\cos \Big(\displaystyle\frac{1}{x^2}\Big)+\nonumber\\
\end{eqnarray*}
$$
+\Big(\displaystyle\frac{24}{x^5}-\underbrace{\displaystyle\frac{36}{x^4}}_{2}
+\underbrace{\displaystyle\frac{12}{x^3}}_{3}\Big)(\ln x)\sin \Big(\displaystyle\frac{1}{x^2}\Big)
-\displaystyle\frac{12}{x^7}
\cos \Big(\displaystyle\frac{1}{x^2}\Big)+\Big(\underbrace{\displaystyle\frac{12}{x^4}}_{4} -\underbrace{\displaystyle\frac{18}{x^5}}_{*}\Big)\sin \Big(\displaystyle\frac{1}{x^2}\Big)+
$$
\begin{eqnarray*}
&&\quad +\left(-\underbrace{\displaystyle\frac{6}{x^5}}_{*}+\underbrace{\displaystyle\frac{12}{x^4}}_{4}\right)\sin \Big(\displaystyle\frac{1}{x^2}\Big)+
\left(\displaystyle\frac{6}{3x^3}-\underbrace{\displaystyle\frac{6}{x^2}}_{5}+\underbrace{\displaystyle\frac{6}{x}}_{6}
\right)\cos \Big(\displaystyle\frac{1}{x^2}\Big)+\nonumber\\
&&\quad +\underbrace{\displaystyle\frac{24}{x^6}(\ln x)
\cos \Big(\displaystyle\frac{1}{x^2}\Big)}_{1}-\Big(\underbrace{\displaystyle\frac{24}{x^3}}_{3}-\underbrace{\displaystyle\frac{36}{x^4}}_{2}\Big)(\ln x)
\sin \Big(\displaystyle\frac{1}{x^2}\Big)+\nonumber
\end{eqnarray*}
$$ 
-\underbrace{\displaystyle\frac{24}{x^4}\sin \Big(\displaystyle\frac{1}{x^2}\Big)}_{4}+
\left(\underbrace{\displaystyle\frac{12}{2x^2}}_{5}-\underbrace{\displaystyle\frac{12}{x}}_{6}\right)
\cos \Big(\displaystyle\frac{1}{x^2}\Big) +\underbrace{\displaystyle\frac{12}{x^3}(\ln x)\sin \Big(\displaystyle\frac{1}{x^2}\Big)}_{3}
+\underbrace{\displaystyle\frac{6}{x}\cos \Big(\displaystyle\frac{1}{x^2}}_{6}\Big)
$$
\begin{eqnarray}\label{eq44}
&=&-\displaystyle\frac{8}{x^9}(\ln x)\sin \Big(\displaystyle\frac{1}{x^2}\Big)+
\displaystyle\frac{36}{x^7}(\ln x)\cos \Big(\displaystyle\frac{1}{x^2}\Big)+
\displaystyle\frac{24}{x^5}(\ln x)\sin \Big(\displaystyle\frac{1}{x^2}\Big)+\nonumber\\
&&\quad -\displaystyle\frac{12}{x^7}\cos \Big(\displaystyle\frac{1}{x^2}\Big)
-\displaystyle\frac{24}{x^5}\sin \Big(\displaystyle\frac{1}{x^2}\Big)+
\displaystyle\frac{2}{x^3}\cos \Big(\displaystyle\frac{1}{x^2}\Big)=f'''(x).
\end{eqnarray}

\medskip
By (\ref{eq42}), (\ref{eq43}) and (\ref{eq44}), we get the first derivative $f'(x)$, the second derivative $f''(x)$ and the third derivative $f'''(x)$ of 
$f(x)=(\ln x)\,\cos\Big(\displaystyle\frac{1}{x^2}\Big)$ by doing 
$\mathcal{R}^{(4)}_{1, 1, 1}$-automatic differentiation and without using any differentiation rule in calculus.

\bigskip
The way of developing the $3$-rd order automatic differentiation induced by   the  $4$-dimensional truncated polynomial real algebra $\mathcal{R}^{(4)}$  in this paper
 can be generalized to establish the mathematical foundation of the $n$-th order automatic differentiation induced by   the  $(n+1)$-dimensional truncated polynomial real algebra $\mathcal{R}^{(n+1)}$ solidly. We will do it in another paper.

\end{document}